\documentclass{amsart}
\usepackage{epsfig}
\newtheorem{theorem}{Theorem}

\newtheorem{lemma}[theorem]{Lemma}
\newtheorem{corollary}[theorem]{Corollary}
\newcommand\beq{\begin{equation}}
\newcommand\eeq{\end{equation}}
\newcommand\bce{\begin{center}}
\newcommand\ece{\end{center}}
\newcommand\bea{\begin{eqnarray}}
\newcommand\eea{\end{eqnarray}}
\newcommand\ben{\begin{enumerate}}
\newcommand\een{\end{enumerate}}

\newcommand\us{\underset}
\newcommand\wt{\widetilde}

\newcommand\nn{\nonumber}

\newcommand\ms{\medskip}
\newcommand\wh{\widehat}
\newcommand\brr{\begin{array}}
\newcommand\err{\end{array}}
\newcommand\bt{\begin{tabular}}
\newcommand\et{\end{tabular}}
\renewcommand\S{{\mathcal S}}
\newcommand\D{{\mathcal D}}
\newcommand\rs{\Psi}
\newcommand\kra{\varphi}
\newcommand\bij{\Phi}
\newcommand\ct{\mathrm{ct}}
\newcommand\lt{\mathrm{lt}}
\newcommand\rt{\mathrm{rt}}
\newcommand\odr{\mathrm{or}}
\newcommand\er{\mathrm{er}}
\newcommand\cmt{\mathrm{cmt}}
\newcommand\ret{\mathrm{ret}}
\newcommand\len{\mathrm{length}}
\newcommand\fp{\mathrm{fp}}
\newcommand\exc{\mathrm{exc}}
\newcommand\des{\mathrm{des}}
\newcommand\oc{\lambda}
\newcommand\C{{\mathbf C}}
\renewcommand\u{{\mathbf u}}
\renewcommand\d{{\mathbf d}}
\newcommand\ih{\mathrm{ih}}
\newcommand\fh{\mathrm{fh}}

\newenvironment{abstrac}{%
         \small
        \begin{center}%
          {\bfseries {Abstract}\vspace{-.5em}}%
                 \end{center}%
        \quotation}

\title{A simple and unusual bijection for Dyck paths and its consequences}

\author{Sergi Elizalde}
\address{Department of Mathematics, MIT, Cambridge, MA 02139.}
\email{sergi@math.mit.edu}
\author{Emeric Deutsch}
\address{Polytechnic University, Brooklyn, NY 11201.}
\email{deutsch@duke.poly.edu}

\begin{document}
\maketitle

\begin{abstrac}
\ms 
In this paper we introduce a new bijection from the set of
Dyck paths to itself. This bijection has the property that it maps
statistics that appeared recently in the study of pattern-avoiding
permutations into classical statistics on Dyck paths, whose
distribution is easy to obtain.

We also present a generalization of the bijection, as well as
several applications of it to enumeration problems of statistics
in restricted permutations.
\end{abstrac}

%%%%%%%%%%%%%%%%%%%%%%%%%%%%%%%%%%%%
\section{Introduction}
%%%%%%%%%%%%%%%%%%%%%%%%%%%%%%%%%%%%

Motivated by the study of statistics on pattern-avoiding
permutations, new statistics on Dyck paths have recently been
introduced in \cite{Eli02,EliPak}. These statistics, which are
based on the notion of tunnel of a Dyck path, have important
applications in the study of pattern-avoiding permutations for
patterns of length 3. Several enumeration problems regarding
permutation statistics can be solved more easily considering their
counterpart in terms of Dyck paths.

In Section~\ref{sec:bij} we present a new bijection $\bij$ from
the set of Dyck paths to itself, and in
Section~\ref{sec:properties} we study its properties. The
interesting properties are that the statistics mentioned above
involving tunnels are mapped by $\bij$ into other known statistics
such as hills, returns and even rises, which have been widely
studied in the literature. This bijection gives new ways to derive
generating functions enumerating those statistics.

In Section~\ref{sec:gener} we give a family of bijections
depending on an integer parameter $r$, from which the main
bijection $\bij$ is the particular case $r=0$. These bijections
give correspondences involving new statistics of Dyck paths, which
generalize the above ones. We give multivariate generating
functions for them.

Finally, Section~\ref{sec:gener} discusses several applications of
our bijections to enumeration of statistics in 321- and
132-avoiding permutations. In particular, we generalize a recent
theorem about fixed points in restricted permutations, and we find
a multivariate generating function for fixed points, excedances
and descents in 132-avoiding permutations.

%%%%%%%%%%%%%%%%%%%%%%%%%%%%%%%%%%%%
\section{Preliminaries}
%%%%%%%%%%%%%%%%%%%%%%%%%%%%%%%%%%%%

Recall that a \emph{Dyck path} of length $2n$ is a lattice path in
$\mathbb{Z}^2$ between $(0,0)$ and $(2n,0)$ consisting of up-steps
$(1,1)$ and down-steps $(1,-1)$ which never goes below the
$x$-axis. Sometimes it will be convenient to encode each up-step
by a letter $\u$ and each down-step by $\d$, obtaining an encoding
of the Dyck path as a \emph{Dyck word}. We shall denote by $\D_n$
the set of Dyck paths of length $2n$, and by
$\D=\bigcup_{n\geq0}\D_n$ the class of all Dyck paths. It is
well-known that $|\D_n|=\C_n=\frac{1}{n+1}{2n \choose n}$, the
$n$-th Catalan number. If $D\in\D_n$, we will write $|D|=n$ to
indicate the semilength of $D$. The generating function that
enumerates Dyck paths according to their semilength is
$\sum_{D\in\D}{z^{|D|}}=\sum_{n\geq0}{\C_n z^n}
=\frac{1-\sqrt{1-4z}}{2z}$, which we denote by $\C(z)$.

We will use $D$ to refer indistinctively to the Dyck path $D$ or
to the Dyck word associated to it. In particular, given
$D_1\in\D_{n_1}$, $\D_2\in\D_{n_2}$, we will write $D_1 D_2$ to
denote the concatenation of $D_1$ and $D_2$ (note that, as seen in
terms of lattice paths, $D_2$ has to be shifted $2n_1$ units to
the right).

A \emph{peak} of a Dyck path $D\in\D$ is an up-step followed by a
down-step (i.e., an occurrence of $\u\d$ in the associated Dyck
word). The $x$-coordinate of a peak is given by the point at the
top of it. A \emph{hill} is a peak at height 1, where the height
is the $y$-coordinate of the top of the peak. Denote by $h(D)$ the
number of hills of $D$. A \emph{valley} of a Dyck path $D\in\D$ is
a down-step followed by an up-step (i.e., an occurrence of $\d\u$
in the associated Dyck word). An \emph{odd rise} is an up-step in
an odd position when the steps are numbered from left to right
starting with 1 (or, equivalently, it is an up-step at odd level
when the steps leaving the x-axis are considered to be at level
1). Denote by $\odr(D)$ the number of odd rises of $D$. \emph{Even
rises} and $\er(D)$ are defined analogously. The $x$-coordinate of
an odd or even rise is given by the rightmost end of the
corresponding up-step.

A \emph{return} of a Dyck path is a down-step landing on the
$x$-axis. An \emph{arch} is a part of the path joining two
consecutive points on the $x$-axis. Clearly for any $D\in\D_n$ the
number of returns equals the number of arches. Denote it by $\ret(D)$.
Define the $x$-coordinate of an arch as the $x$-coordinate of its leftmost point.

We will use the concept of tunnel introduced by the first author
in~\cite{Eli02}. For any $D\in\D$, define a \emph{tunnel} of $D$
to be a horizontal segment between two lattice points of $D$ that
intersects $D$ only in these two points, and stays always below
$D$. Tunnels are in obvious one-to-one correspondence with
decompositions of the Dyck word $D=A\u B\d C$, where $B\in\D$ (no
restrictions on $A$ and $C$). In the decomposition, the tunnel is
the segment that goes from the beginning of the $\u $ to the end
of the $\d$. If $D\in\D_n$, then $D$ has exactly $n$ tunnels,
since such a decomposition can be given for each up-step $\u $ of
$D$.

A tunnel of $D\in\D_n$ is called a \emph{centered tunnel} if the
$x$-coordinate of its midpoint is $n$, that is, the tunnel is
centered with respect to the vertical line through the middle of
$D$. In terms of the decomposition $D=A\u B\d C$, this is
equivalent to saying that $A$ and $C$ have the same length. Denote
by $\ct(D)$ the number of centered tunnels of $D$.

\begin{figure}[hbt]
\epsfig{file=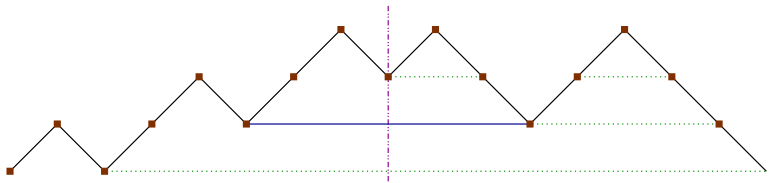,width=4in} \caption{\label{fig:ctrt} One
centered and four right tunnels.}
\end{figure}

A tunnel of $D\in\D_n$ is called a \emph{right tunnel} if the
$x$-coordinate of its midpoint is strictly greater than $n$, that
is, the midpoint of the tunnel is to the right of the vertical
line through the middle of $D$. Clearly, in terms of the
decomposition $D=A\u B\d C$, this is equivalent to saying that the
length of $A$ is strictly bigger than the length of $C$. Denote by
$\rt(D)$ the number of right tunnels of $D$. In
Figure~\ref{fig:ctrt}, there is one centered tunnel drawn with a
solid line, and four right tunnels drawn with dotted lines.
Similarly, a tunnel is called a \emph{left tunnel} if the
$x$-coordinate of its midpoint is strictly less than $n$. Denote
by $\lt(D)$ the number of left tunnels of $D$. Clearly,
$\lt(D)+\rt(D)+\ct(D)=n$ for any~$D\in\D_n$.

For any $D\in\D$, we define a \emph{multitunnel} of $D$ to be a
horizontal segment between two lattice points of $D$ such that $D$
never goes below it. In other words, a multitunnel is just a
concatenation of tunnels, so that each tunnel starts at the point
where the previous one ends. Similarly to the case of tunnels,
multitunnels are in obvious one-to-one correspondence with
decompositions of the Dyck word $D=ABC$, where $B\in\D$ is not
empty. In the decomposition, the multitunnel is the segment that
connects the initial and final points of $B$.

A multitunnel of $D\in\D_n$ is called a \emph{centered
multitunnel} if the $x$-coordinate of its midpoint (as a segment)
is $n$, that is, the tunnel is centered with respect to the
vertical line through the middle of $D$. In terms of the
decomposition $D=ABC$, this is equivalent to saying that $A$ and
$C$ have the same length. Denote by $\cmt(D)$ the number of
centered multitunnels of $D$.

\begin{figure}[hbt]
\epsfig{file=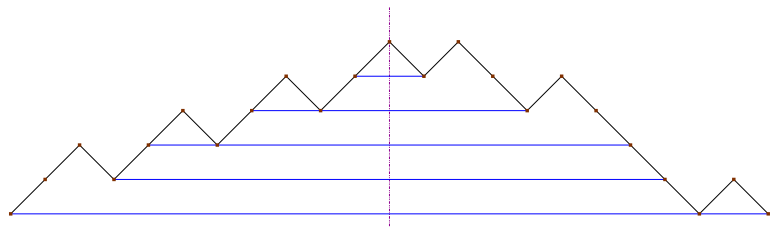,width=4in} \caption{\label{fig:mult} Five
centered multitunnels, two of which are centered tunnels.}
\end{figure}

%%%%%%%%%%%%%%%%%%%%%%%%%%%%%%%
\section{The bijection}\label{sec:bij}
%%%%%%%%%%%%%%%%%%%%%%%%%%%%%%%

In this section we describe a bijection $\bij$ from $\D_n$ to
itself. Let $D\in\D_n$. Each up-step of $D$ has a corresponding
down-step together with which it determines a tunnel. Match each
such pair of steps. Let $\sigma\in\S_{2n}$ be the permutation
defined by
$$\sigma_i=\left\{ \brr{cc} {i+1 \over 2} &\mbox{if $i$ is odd;}
\\ 2n+1-{i\over 2} &\mbox{if $i$ is even.} \err \right. $$
In two-line notation,

$$ \sigma = \left( \brr{ccccccccccc}
1 & 2 & 3 & 4 & 5 & 6 & \cdots & 2n-3 & 2n-2 & 2n-1 & 2n \\
1 & 2n & 2 & 2n-1 & 3 & 2n-2 & \cdots & n-1 & n+2 & n & n+1 \err \right).$$

Then $\bij(D)$ is created as follows. For $i$ from 1 to $2n$,
consider the $\sigma_i$-th step of $D$ (i.e., $D$ is read in zigzag). If
its corresponding matching step has not yet been read, define the
$i$-th step of $\bij(D)$ to be an up-step, otherwise let it be a
down-step. In the first case, we say that the $\sigma_i$-th step
of $D$ \emph{opens} a tunnel, in the second we say that it
\emph{closes} a tunnel.

\begin{figure}[hbt]
\epsfig{file=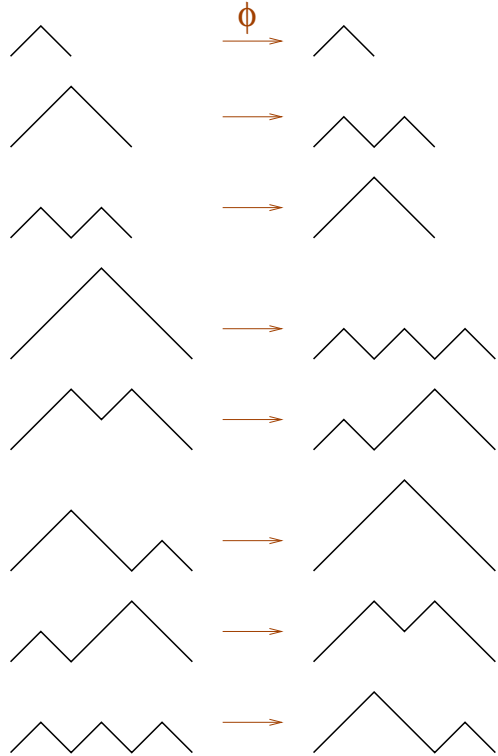,height=4in} \caption{\label{fig:bij125} The bijection $\bij$ for paths of length at most 3.}
\end{figure}

The bijection $\bij$ applied to the Dyck
paths of semilength at most 3 is shown in Figure~\ref{fig:bij125}.

\begin{figure}[hbt]
\epsfig{file=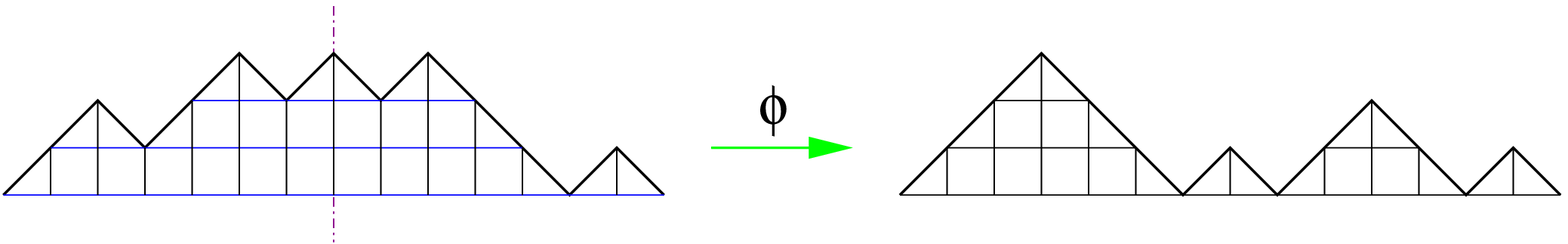,width=6in} \caption{\label{fig:bij} An example of $\bij$.}
\end{figure}

Figure~\ref{fig:bij} shows $\bij$ applied to the example of the Dyck path $D=\u
\u \d \u \u \d \u \d \u \d \d \d \u \d $.

It is clear from the definition that $\bij(D)$ is a Dyck path.
Indeed, it never goes below the x-axis because at any point the
number of down-steps drawn so far can never exceed the number of
up-steps, since each down-step is drawn when the second step of a
matching pair in $D$ is read, and in that case the first step of
the pair has already produced an up-step in $\bij(D)$. Also,
$\bij(D)$ ends in $(2n,0)$ because each of the matched pairs of
$D$ produces an up-step and a down-step in $\bij(D)$.

\begin{figure}[hbt]
\epsfig{file=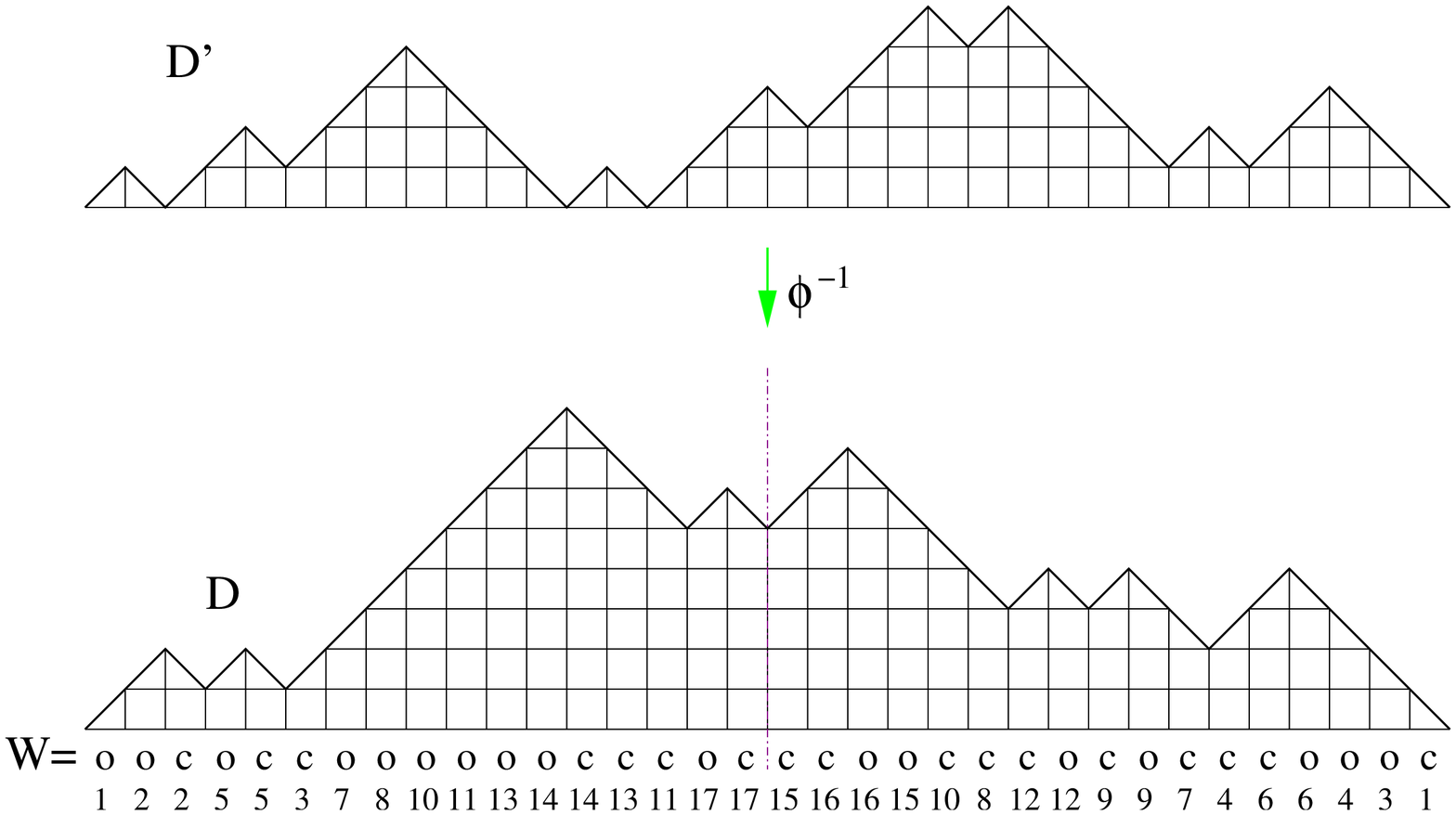,width=5.5in} \caption{\label{fig:inv} The inverse of $\bij$.}
\end{figure}

To show that $\bij$ is indeed a bijection, we will describe the
inverse map $\bij^{-1}$. Given $D'\in\D_n$, the following
procedure recovers the $D\in\D_n$ such that $\bij(D)=D'$. Consider
the permutation $\sigma$ defined above, and let $W=w_1 w_2 \cdots
w_{2n}$ be the word obtained from $D'$ as follows. For $i$ from 1
to $2n$, if the $i$-th step of $D'$ is an up-step, let
$w_{\sigma_i}=o$, otherwise let $w_{\sigma_i}=c$. $W$ contains the
same information as $D'$, with the advantage that the $o$'s are
located in the positions of $D$ in which a tunnel is opened when
$D$ is read in zigzag, and the $c$'s are located in the positions
where a tunnel is closed. Equivalently, the $o$'s are located in
the positions of the left walls of the left and centered tunnels
of $D$, and in the positions of the right walls of the right
tunnels. For an example see Figure~\ref{fig:inv}.

Now we define a matching between the $o$'s and the $c$'s in $W$,
so that each matched pair will give a tunnel in $D$. We will label
the $o$'s with $1,2,\ldots,n$ and similarly the $c$'s, to indicate
that an $o$ and a $c$ with the same label are matched. By left
(resp. right) half of $W$ we mean the symbols $w_i$ with $i\leq n$
(resp. $i>n$). For $i$ from 1 to $2n$, if $w_{\sigma_i}=o$, place
in it the smallest label that has not been used yet. If
$w_{\sigma_i}=c$, match it with the unmatched $o$ in the same half
of $W$ as $w_{\sigma_i}$ with largest label, if such an $o$
exists. If it does not, match $w_{\sigma_i}$ with the unmatched
$o$ in the opposite half of $W$ with smallest label. Note that
since $D'$ was a Dyck path, at any time the number of $c$'s read
so far does not exceed the number of $o$'s, so each $c$ has some
$o$ to be paired up with.

Once the symbols in $W$ have been labelled, $D$ can be recovered
by reading the labels from left to right, drawing an up-step for
each label that is read for the first time, and a down-step for
each label that appears the second time. In Figure~\ref{fig:inv} the labelling is shown under $W$.

%%%%%%%%%%%%%%%%%%%%%%%%%%%%%%%%%%
\section{Properties of $\bij$}\label{sec:properties}
%%%%%%%%%%%%%%%%%%%%%%%%%%%%%%%%%%

\begin{lemma} \label{lemma} Let $D=ABC$ be a decomposition of a Dyck path $D$,
where $B$ is a Dyck path, and $A$ and $C$ have the same length.
Then $\bij(ABC)=\bij(AC)\bij(B)$. In particular, $\bij(\u B\d )=\u
\d \bij(B)$.
\end{lemma}

\begin{proof} It follows immediately from the definition of
$\bij$, since the path $D$ is read in zigzag while $\bij(D)$ is
built from left to right.
\end{proof}

\begin{theorem} \label{corresp} Let $D$ be any Dyck path, and let $D'=\bij(D)$. We have the following
correspondences:
\ben
\item $\ct(D)=h(D')$,
\item $\rt(D)=\er(D')$,
\item $\lt(D)+\ct(D)=\odr(D')$,
\item $\cmt(D)=\ret(D')$.
\een\end{theorem}

\begin{proof}
First we show (1). Consider a centered tunnel given by the
decomposition $D=A\u B\d C$. Applying Lemma~\ref{lemma} twice,
$D'=\bij(A\u B\d C)=\bij(AC)\bij(\u B\d )=\bij(AC)\u \d \bij(B)$,
so we have a hill $\u \d $ in $D'$. Reciprocally, any hill in
$D'$, say $D'=X\u \d Y$, where $X,Y\in\D$, comes from a centered
tunnel $D=Z_1 \u  \bij^{-1}(Y) \d  Z_2$, where $Z_1$ and $Z_2$ are
respectively the first and second halves of $\bij^{-1}(X)$.

The proof of (4) is very similar. Recall that $\ret(D')$ equals
the number of arches of $D'$. Given a centered multitunnel
corresponding to the decomposition $D=ABC$, we have
$\bij(D)=\bij(AC)\bij(B)$, so $D'$ has an arch starting at the
first step of $\bij(B)$, which is nonempty.

To show (2), consider a right tunnel given by the decomposition
$D=A\u B\d C$, where $\len(A)>\len(C)$. Of the two steps $\u $ and
$\d $ delimiting the tunnel, $\d $ will be encountered before $\u
$ when $D$ is read in zigzag, since $\len(A)>\len(C)$. So $\d $
will open a tunnel, producing an up-step in $D'$. Besides, this
up-step will be at an even position, since $\d $ was in the right
half of $D$. Reciprocally, an even rise of $D'$ corresponds to a
step in the right half of $D$ that opens a tunnel when $D$ is read
in zigzag, so it is necessarily a right tunnel.

Relation (3) follows from (2) and the fact that the total number
of tunnels of $D$ is $\lt(D)+\ct(D)+\rt(D)=n$, and the total
number of up-steps of $D'$ is $\odr(D')+\er(D')=n$.
\end{proof}

\ms

One of the most interesting applications of this bijection is that
it can be used to enumerate Dyck paths according to the number of
centered, left, and right tunnels, and number of centered
multitunnels. We are looking for a multivariate generating
function for these four statistics, namely
$$\wt{G}(x,u,v,w,z)=\sum_{D\in\D} x^{\ct(D)} u^{\lt(D)} v^{\rt(D)} w^{\cmt(D)} z^{|D|}.$$
By Theorem~\ref{corresp}, this generating function can be
expressed as $\wt{G}(x,u,v,w,z)=G(\frac{x}{u},u,v,w,z)$, where
$$G(t,u,v,w,z)=\sum_{D\in\D} t^{h(D)} u^{\odr(D)} v^{\er(D)} w^{\ret(D)} z^{|D|}.$$

We can derive an equation for $G$ using the symbolic method
described in \cite{FlSe98} and \cite{SeFl96}. A recursive
definition for the class $\D$ is given by the fact that every
nonempty Dyck path $D$ can be decomposed in a unique way as $D=\u
A\d B$, where $A,B\in\D$. The number of hills of $\u A\d B$ is
$h(B)+1$ if $A$ is empty, and $h(B)$ otherwise. The odd rises of
$A$ become even rises of $\u A\d B$, and the even rises of $A$
become odd rises of $\u A\d B$. Thus, we have $\er(\u A\d
B)=\odr(A)+\er(B)$, and $\odr(\u A\d B)=\er(A)+\odr(B)+1$, where
the extra odd rise comes from the first step $\u $. We also have
$\ret(\u A\d B)=\ret(B)+1$. Hence, we obtain the following
equation for $G$: \bea
G(t,u,v,w,z)=1+uzw(G(1,v,u,1,z)-1+t)G(t,u,v,w,z). \label{eqG} \eea
Denote $G_1:=G(1,u,v,1,z)$, $H_1:=G(1,v,u,1,z)$. Substituting
$t=w=1$ in (\ref{eqG}), we obtain \bea G_1=1+uz H_1 G_1,
\label{eqG1} \eea and interchanging $u$ and $v$, \bea H_1=1+vz G_1
H_1. \label{eqH1} \eea
 Solving (\ref{eqG1}) and (\ref{eqH1}) for $H_1$, gives
$$H_1=\frac{1+(u-v)z-\sqrt{1-2(v+u)z+(v-u)^2 z^2}}{2uz}.$$
Thus, from (\ref{eqG}), \bea
G(t,u,v,w,z)&=&\frac{1}{1-uzw(H_1-1+t)}\nn \\
&=&\frac{2}{2-w+(v+u-2tu)wz+w\sqrt{1-2(v+u)z+(v-u)^2 z^2}}.
\nn\eea

Now, switching to $\wt{G}$, we obtain the following theorem.
\begin{theorem} \label{gftun} The multivariate generating function for Dyck paths according to centered, left, and right tunnels, centered multitunnels, and semilength is
$$\sum_{D\in\D} x^{\ct(D)} u^{\lt(D)} v^{\rt(D)} w^{\cmt(D)} z^{|D|}=\frac{2}{2-w+(v+u-2x)wz+w\sqrt{1-2(v+u)z+(v-u)^2 z^2}}.$$
\end{theorem}

%%%%%%%%%%%%%%%%%%%%%%%%%%%%%%%%%%%%
\section{Generalizations}\label{sec:gener}
%%%%%%%%%%%%%%%%%%%%%%%%%%%%%%%%%%%%

Here we present a generalization $\bij_r$ of the bijection $\bij$,
which depends on a nonnegative integer parameter $r$. Given
$D\in\D_n$, copy the first $2r$ steps of $D$ into the first $2r$
steps of $\bij_r(D)$. Now, read the remaining steps of $D$ in
zigzag in the following order: $2r+1$, $2n$, $2r+2$, $2n-1$,
$2r+3$, $2n-2$, and so on. For each of these steps, if its
corresponding matching step in $D$ has not yet been encountered,
draw an up-step in $\bij_r(D)$, otherwise draw a down-step. Note
that for $r=0$ we get the same bijection $\bij$ as before.

Note that $\bij_r$ can be defined exactly as $\bij$ with the
difference that instead of $\sigma$, the permutation that gives
the order in which the steps of $D$ are read is
$\sigma^{(r)}\in\S_{2n}$, defined as
$$\sigma^{(r)}_i=\left\{ \brr{cc} i &\mbox{ if $i\le 2r$;} \\ {i+1 \over 2}+r &\mbox{ if $i>2r$ and $i$ is odd;}
\\ 2n+1-{i\over 2}-r &\mbox{ if $i>2r$ and $i$ is even.} \err \right. $$

\begin{figure}[hbt]
\epsfig{file=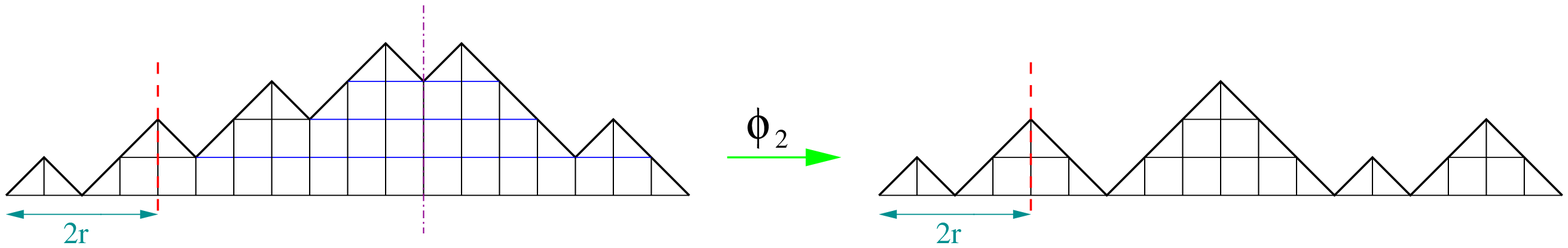,width=6in} \caption{\label{fig:bij_r} An example of $\bij_2$.}
\end{figure}

Figure~\ref{fig:bij_r} shows an example of $\bij_r$ for $r=2$ applied to the path $D=\u\d\u\u\d\u\u\d\u\u\d\u\d\d\d\u\d\d$.

It is clear from the definition that $\bij_r(D)$ is a Dyck path. A
reasoning similar to the one used for $\bij$ shows that $\bij_r$
is indeed a bijection.

The properties of $\bij$ given in Theorem~\ref{corresp} generalize
to analogous properties of $\bij_r$. We will prove them using the
following lemma, which follows immediately from the definition of
$\bij_r$.

\begin{lemma} \label{lemma_r} Let $r\geq 0$, and let $D=ABC$ be a decomposition of a Dyck path $D$,
where $B$ is a Dyck path, and $\len(A)=\len(C)+2r$. Then
$\bij_r(ABC)=\bij_r(AC)\bij(B)$.
\end{lemma}

\begin{theorem} \label{corresp_r} Let $r\geq 0$, let $D$ be any Dyck path, and let $D'=\bij_r(D)$. We have the following
correspondences: \ben
\item $\#\{\mbox{tunnels of }D \mbox{ with midpoint at } x=n+r \}\ =\ \#\{\mbox{hills of }
D' \mbox{ in } x>2r\}$,
\item $\#\{\mbox{tunnels of }D \mbox{ with midpoint in } x>n+r \}\ =\ \#\{\mbox{even rises of } D'
\mbox{ in } x>2r\}$,
\item $\#\{\mbox{tunnels of }D \mbox{ with midpoint in } x\le n+r \}\ =\ \#\{\mbox{odd rises of } D'
\mbox{ in } x>2r\}\\ \hspace*{8cm} +\#\{\mbox{up-steps of } D'
\mbox{ in } x\le 2r\}$,
\item $\#\{\mbox{multitunnels of }D \mbox{ with midpoint at } x=n+r\}\ =\ \#\{\mbox{arches of } D'
\mbox{ in } x\ge 2r\}$. \een\end{theorem}

\begin{proof}
Fist we show (1). A tunnel given by the decomposition $D=A\u B\d
C$ has its midpoint at $x=n+r$ exactly when $\len(A)=\len(C)+2r$.
Applying Lemmas \ref{lemma_r} and ~\ref{lemma}, $D'=\bij_r(A\u B\d
C)=\bij_r(AC)\bij(\u B\d )=\bij_r(AC)\bij(\u \d
)\bij(B)=\bij_r(AC)\u \d \bij(B)$, and $\u \d $ is a hill of $D'$
in $x>2r$, since $\len(\bij_r(AC))\geq 2r$. Reciprocally, any hill
of $D'$ in $x>2r$, say $D'=X\u \d Y$, where $X,Y\in\D$ and
$\len(X)\geq 2r$, comes from a tunnel with midpoint at $x=n+r$,
namely $D=Z_1 \u \bij^{-1}(Y) \d  Z_2$, where $Z_1
Z_2=\bij_r^{-1}(X)$ and $\len(Z_1)=\len(Z_2)+2r$.

The proof of (4) is very similar. A multitunnel given by $D=ABC$
has its midpoint at $x=n+r$ exactly when $\len(A)=\len(C)+2r$. In
this case, $\bij_r(D)=\bij_r(AC)\bij(B)$ by Lemma~\ref{lemma_r},
so $D'$ has an arch starting at the first step of $\bij(B)$. Notice that this arch is
in $x\geq 2r$ because $\len(\bij_r(AC))\geq 2r$.

To show (2), consider a tunnel in $D$ with midpoint in $x>n+r$.
This is equivalent to saying that it is given by a decomposition
$D=A\u B\d C$ with $\len(A)>\len(C)+2r$. In particular, the tunnel
is contained in the halfspace $x\geq 2r$, so the two steps $\u $
and $\d $ delimiting the tunnel are in the part of $D$ that is
read in zigzag in the process to obtain $\bij_r(D)$, and $\d $
will be encountered before $\u $, since $\len(A)-2r>\len(C)$. So
$\d $ will open a tunnel, producing an up-step of $D'$ in $x>2r$.
Besides, this up-step will be at an even position, since $\d $ is
in $x>n+r$, that is, in the right half of the part of $D$ that is
read in zigzag. Reciprocally, an even rise of $D'$ in $x>2r$
corresponds to a step of $D$ in $x>n+r$ that opens a tunnel when
$D$ is read according to $\sigma^{(r)}$, so it is necessarily
tunnel with midpoint to the right of $x=n+r$ .

Relation (3) follows from (2) and the fact that the total number
of tunnels of $D$ is \\ $\#\{\mbox{tunnels of }D \mbox{ with
midpoint in } x>n+r \}+\#\{\mbox{tunnels of }D \mbox{ with
midpoint in } x\le n+r \}=n$, and the total number of up-steps of
$D'$ is $\#\{\mbox{even rises of } D' \mbox{ in }
x>2r\}+\#\{\mbox{odd rises of } D' \mbox{ in }
x>2r\}+\#\{\mbox{up-steps of } D' \mbox{ in } x\le 2r\}=n$.
\end{proof}
%\ms

Similarly to how we used the properties of $\bij$ to prove
Theorem~\ref{gftun}, we can use the properties of $\bij_r$ to
prove a more general theorem. Our goal is to enumerate Dyck paths
according to the number tunnels with midpoint on, to the right of,
and to the left of an arbitrary vertical line $x=n+r$, and
multitunnels with midpoint on that line. In generating function
terms, we are looking for an
expression for \bea\nn F(t,u,v,w,y,z):=\hspace{11.5cm} \\
\nn=\us{0\le r\le n}{\sum_{n\ge0}} \sum_{D\in\D_n}
t^{\#\{\mbox{tun. of $D$ w/ midp. at $x=n+r$}\}}
u^{\#\{\mbox{tun. of $D$ w/ midp. in $x\le n+r$}\}}\\
\nn v^{\#\{\mbox{tun. of $D$ w/ midp. in $x>n+r$}\}}
w^{\#\{\mbox{multitun. of $D$ w/ midp. at $x=n+r$}\}} y^r z^n.
\eea

Note that the variable $y$ marks the position of the vertical line
$x=n+r$ with respect to which the tunnels are classified. The
following theorem gives an expression for $F$.

\begin{theorem} Let $F$, $G$ and $\C$ be defined as above. Then,
\bea\nn F(t,u,v,w,y,z)=\frac{\C(uyz)G(t,u,v,w,z)}{1-y u^2 z^2
\C^2(uyz) G(1,u,v,1,z) G(1,v,u,1,z)}\\ \nn
=\frac{2B(2+(v-u)z+A)}{\left[2+(u+v-2tu)wz+wA\right]\left[(A+(v-u)z)B-4uyz\right]}
,\eea where $A:=\sqrt{1-2(u+v)z+(u-v)^2 z^2}-1$,
$B:=\sqrt{1-4uyz}-1$.\end{theorem}

\begin{proof}
By Theorem~\ref{corresp_r}, the generating function $F$ can be
expressed as \bea\nn F(t,u,v,w,y,z)=\hspace{11.6cm} \\
\nn\us{0\le r\le n}{\sum_{n\ge0}} \sum_{D\in\D_n}
t^{\#\{\mbox{hills of $D$ in $x>2r$}\}} u^{\#\{\mbox{odd rises of
$D$ in $x>2r$}\}+
\#\{\mbox{up-steps of $D$ in $x\le 2r$}\}}\\
v^{\#\{\mbox{even rises of $D$ in $x>2r$}\}} w^{\#\{\mbox{arches
of $D$ in $x\ge 2r$}\}} y^r z^n. \label{Frises} \eea

For each path $D$ in this summation, the $y$-coordinate of its
intersection with the vertical line $x=2r$ has to be even. Fix
$h\ge 0$. We will now focus only on the paths $D\in\D$ for which
this intersection has $y$-coordinate equal to $2h$. Let $D=AB$,
where $A$ and $B$ are the parts of the path respectively to the
left and to the right of $x=2r$. Then, $\#\{\mbox{hills of } D
\mbox{ in } x>2r\}=\#\{\mbox{hills of } B\}$, $\#\{\mbox{odd rises
of } D \mbox{ in } x>2r\}=\#\{\mbox{odd rises of } B\}$,
$\#\{\mbox{up-steps of } D \mbox{ in } x\le
2r\}=\#\{\mbox{up-steps of } A\}$, and $\#\{\mbox{arches of } D
\mbox{ in } x\ge 2r\}=\#\{\mbox{arches of } B\}$.

$B$ can be any path starting at height $2h$ and landing on the
$x$-axis, never going below it. If $h>0$, consider the first
down-step of $B$ that lands at height $2h-1$. Then $B$ can be
decomposed as $B=B_1dB'$, where $B_1$ is any Dyck path, and $B'$
is any path starting at height $2h-1$ and landing on the $x$-axis,
never going below it. Applying this decomposition recursively, $B$
can be written uniquely as $B=B_1\d B_2\d \cdots B_{2h}\d
B_{2h+1}$, where the $B_i$'s for $1\le i\le 2h+1$ are arbitrary
Dyck paths. The number of hills and number of arches of $B$ are
given by those of $B_{2h+1}$. The odd rises of $B$ are the odd
rises of the $B_i$'s with odd subindex plus the even rises of
those with even subindex. In a similar way one can describe the
even rises of $B$. The semilength of $B$ is the sum of semilengths
of the $B_i$'s plus $h$, which comes from the $2h$ additional
down-steps. Thus, the generating function for all paths $B$ of
this form, where $t$, $u$, $v$, and $z$ mark respectively number
of hills, number of odd rises, number of even rises, and
semilength, is \bea z^h G^h(1,u,v,1,z) G^h(1,v,u,1,z)
G(t,u,v,w,z).\label{pathB}\eea

Similarly, $A$ can be decomposed uniquely as $A=A_1\u A_2\u \cdots
A_{2h}\u A_{2h+1}$. The number of up-steps of $A$ is the sum of
the number of up-steps of each $A_i$, plus a $2h$ term that comes
from the additional up-steps. The generating function for paths
$A$ of this form, where $u$ marks the number of up-steps, and $y$
and $z$ mark both the semilength, is \bea z^h y^h u^{2h}
\C^{2h+1}(uyz) .\label{pathA}\eea

The product of (\ref{pathB}) and (\ref{pathA}) gives the
generating function for paths $D=AB$ where the height of the
intersection point of $D$ with the vertical line between $A$ and
$B$ is $2h$, where the variables mark the same statistics as in
(\ref{Frises}). Note that the exponent of $y$ is half the distance
between the origin of $D$ and this vertical line. Summing over
$h$, we obtain \bea \nn F(t,u,v,w,y,z)=\sum_{h\ge 0} z^{2h}
y^h u^{2h} \C^{2h+1}(uyz) G^h(1,u,v,1,z) G^h(1,v,u,1,z) G(t,u,v,w,z)\\
\nn = \frac{\C(uyz)G(t,u,v,w,z)}{1-y u^2 z^2 \C^2(uyz)
G(1,u,v,1,z) G(1,v,u,1,z)}. \eea
\end{proof}

%%%%%%%%%%%%%%%%%%%%%%%%%%%%%%%%%%%%%%%%%%%%%%%%
\section{Connection to pattern-avoiding
permutations}\label{sec:perm}
%%%%%%%%%%%%%%%%%%%%%%%%%%%%%%%%%%%%%%%%%%%%%%%%

The bijection $\bij$ has applications to the subject of
enumeration of statistics in pattern-avoiding permutations.

First we review some basic definitions in this subject. Given two
permutations $\pi=\pi_1\pi_2\cdots\pi_n\in\S_n$ and
$\sigma=\sigma_1\sigma_2\cdots\sigma_m\in\S_m$, with $m\le n$, we
say that $\pi$ \emph{contains} $\sigma$ if there exist indices
$i_1<i_2<\ldots<i_m$ such that $\pi_{i_1}\pi_{i_2}\cdots\pi_{i_m}$
is in the same relative order as $\sigma_1\sigma_2\cdots\sigma_m$.
If $\pi$ does not contain $\sigma$, we say that $\pi$ is
\emph{$\sigma$-avoiding}. For example, if $\sigma=132$, then
$\pi=24531$ contains $\sigma$, because $\pi_1\pi_3\pi_4=253$.
However, $\pi=42351$ is $\sigma$-avoiding.

We say that $i$ is a \emph{fixed point} of a permutation $\pi$ if
$\pi_i=i$, and that $i$ is an \emph{excedance} of $\pi$ if
$\pi_i>i$. We say that $i\le n-1$ is a \emph{descent} of
$\pi\in\S_n$ if $\pi_i>\pi_{i+1}$. Denote by $\fp(\pi)$,
$\exc(\pi)$, and $\des(\pi)$ the number of fixed points, the
number of excedances, and the number of descents of $\pi$
respectively. Denote by $\S_n(\sigma)$ the set of
$\sigma$-avoiding permutations in $\S_n$.

The distribution of fixed points in pattern-avoiding permutations
was considered for the first time in \cite{RSZ02}. There appears
the following result.

\begin{theorem}[\cite{RSZ02}] \label{th:RSZ}
The number of 321-avoiding permutations $\pi \in \S_n$ with
$\fp(\pi)=k$ equals the number of 132-avoiding
permutations $\pi \in \S_n$ with $\fp(\pi)=k$,
for any  $0 \le k\le n$.
\end{theorem}

This theorem can be refined considering not only fixed points but
also excedances, as shown in \cite{Eli02}.

\begin{theorem}[\cite{Eli02}] \label{th:old}
The number of 321-avoiding permutations $\pi \in \S_n$ with
$\fp(\pi)=k$ and $\exc(\pi)=l$ equals the number of 132-avoiding
permutations $\pi \in \S_n$ with $\fp(\pi)=k$ and $\exc(\pi)=l$,
for any  $0 \le k,l\le n$.
\end{theorem}

The proofs given in \cite{RSZ02,Eli02} are not bijective. The
first bijective proof of these theorems appears in \cite{EliPak}.
Our bijection $\bij_r$ can be used to give a bijective proof of
the following generalization of Theorem~\ref{th:RSZ}. Note that
the particular case $r=0$ gives a new bijective proof of such
theorem. 

\begin{theorem} Fix $r,n\ge 0$. For any $\pi\in\S_n$, define
$\alpha_r(\pi)=\#\{i:\pi_i=i+r\}$,
$\beta_r(\pi)=\#\{i:i>r,\pi_i=i\}$. Then, the number of
321-avoiding permutations $\pi \in \S_n$ with $\beta_r(\pi)=k$
equals the number of 132-avoiding permutations $\pi \in \S_n$ with
$\alpha_r(\pi)=k$, for any  $0 \le k\le n$.
\end{theorem}

\begin{proof}
We will use two bijections given in \cite{Eli02}, one from
321-avoiding permutations to Dyck paths and one from 132-avoiding
permutations to Dyck paths.

The first bijection, which we denote
$\rs:\S_n(321)\longrightarrow\D_n$, can be defined as follows. Any
permutation $\pi\in\S_n$ can be represented as an $n\times n$
array with a cross on the squares $(i,\pi_i)$. It is known that
$\pi$ is 321-avoiding if and only if both the subsequence
determined by its excedances and the one determined by the
remaining elements are increasing. Given the array of
$\pi\in\S_n(321)$, consider the path with \emph{down} and
\emph{right} steps along the edges of the squares that goes from
the upper-left corner to the lower-right corner of the array
leaving all the crosses to the right and remaining always as close
to the main diagonal as possible. Define $\rs(\pi)$ to be the Dyck
path obtained from this path by reading an up-step every time the
path moves down, and a down-step every time the path moves to the
right. Figure~\ref{fig:bij_rs} shows a picture of this bijection
for $\pi=23147586$.

\begin{figure}[hbt]
\epsfig{file=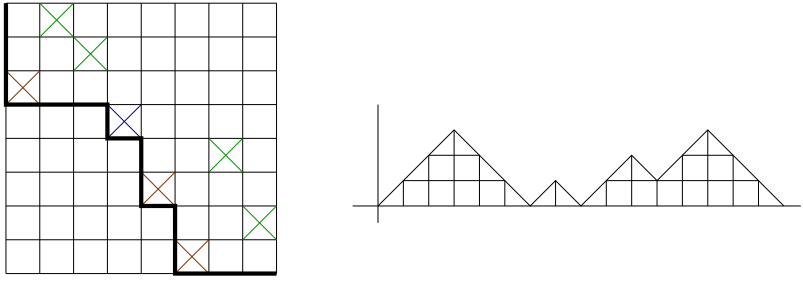,height=3.7cm} \caption{\label{fig:bij_rs}
The bijection $\rs$.}
\end{figure}

It can easily be checked that $\rs$ maps fixed points to hills.
More precisely, $i$ is a fixed point of $\pi$ if and only if
$\rs(\pi)$ has a hill with $x$-coordinate $2i-1$. This implies
that $\beta_r(\pi)=\#\{\mbox{hills of } \rs(\pi) \mbox{ in }
x>2r\}$.

The second bijection, which we denote
$\kra:\S_n(132)\longrightarrow\D_n$, was given by Krattenthaler in
\cite{Kra01}, up to reflection of the path over a vertical line.
It can be described as follows. Given $\pi\in\S_n(132)$
represented by an $n\times n$ array as before, consider the path
with \emph{up} and \emph{right} steps that goes from the
lower-left corner to the upper-right corner, leaving all the
crosses to the right, and staying always as close to the diagonal
connecting these two corners as possible. Then $\kra(\pi)$ is the
Dyck path obtained from this path by reading an up-step every time
the path goes up and a down-step every time it goes right.
Figure~\ref{fig:bij2} shows an example when $\pi=67435281$.

\begin{figure}[hbt]
\epsfig{file=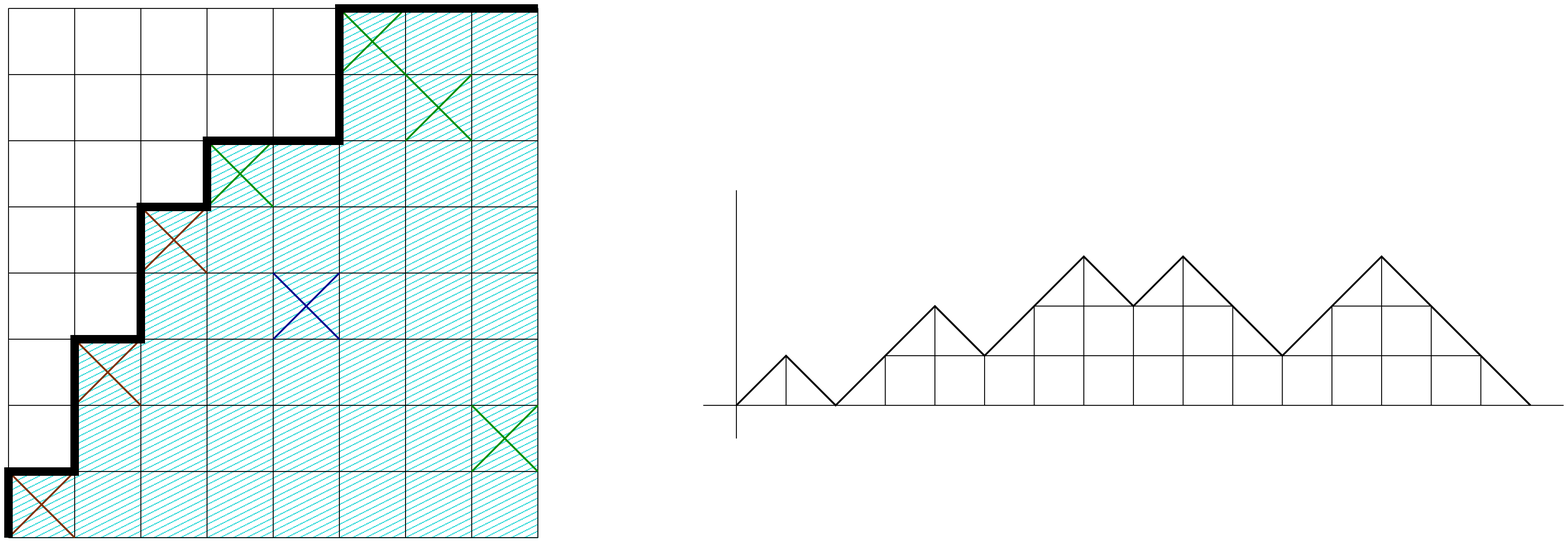,height=3.7cm} \caption{\label{fig:bij2}
The bijection $\kra$.}
\end{figure}

The property of $\kra$ that is useful here is that an element $i$
with $\pi_i=i+r$ corresponds to a tunnel of $\kra(\pi)$ with
midpoint at $x=n+r$. To show this, we repeat the reasoning in
\cite{Eli02,EliPak}. There is an easy way to recover a permutation
$\pi\in\S_n(132)$ from $\kra(\pi)$. Consider the path from the
lower-left corner to the upper-right corner or the array, used to
define $\kra(\pi)$. Now, row by row, put a cross in the leftmost
square to the right of this path such that there is exactly one
cross in each column. This gives us $\pi$ back. To each cross we
can associate a tunnel in a natural way. Indeed, each cross
produces a decomposition $\kra(\sigma)=A\u B\d C$ where $B$
corresponds to the part of the path above and to the left of the
cross. Here $\u $ corresponds to the vertical step directly to the
left of the cross, and $\d$ to the horizontal step directly above
the cross. Thus, fixed points, which correspond to crosses on the
main diagonal, give centered tunnels. More generally, crosses
$(i,i+r)$ give tunnels with midpoint $r$ units to the right of the
center, that is, tunnels of $\kra(\pi)$ with midpoint at $x=n+r$.
Therefore, we have that $\alpha_r(\pi)=\#\{\mbox{tunnels of
}\kra(\pi) \mbox{ with midpoint at } x=n+r \}$.

Now all we need to do is use $\bij_r$ and property (1) given in
Theorem~\ref{corresp_r}. From this it follows that the bijection
$\rs^{-1}\circ\bij_r\circ\kra:\S_n(132)\longrightarrow\S_n(321)$
has the property that
$\beta_r(\rs^{-1}\circ\bij_r\circ\kra(\pi))=\#\{\mbox{hills of }
\bij_r\circ\kra(\pi) \mbox{ in } x>2r\}=\#\{\mbox{tunnels of
}\kra(\pi) \mbox{ with midpoint at } x=n+r \}=\alpha_r(\pi)$.
\end{proof}

\ms

The proof of Theorem~\ref{th:old} given in \cite{Eli02} describes
a simple way to enumerate 321-avoiding permutations with respect
to the statistics $\fp$ and $\exc$. However, the analogous
enumeration for 132-avoiding permutations is done in a more
intricate way. As an application of the properties of $\bij$, we
give a more direct derivation of the multivariate generating
function for 132-avoiding permutations according to number of
fixed points and number of excedances.

\begin{corollary}[of Theorem~\ref{gftun}]\label{cor}
\bea\label{fpexc132} \sum_{n\ge 0}
\sum_{\pi\in\S_n(132)}{x^{\fp(\pi)}v^{\exc(\pi)}}z^n=\frac{2}{1+(1+v-2x)z+\sqrt{1-2(1+v)z+(1-v)^2
z^2}}.\eea
\end{corollary}

\begin{proof}
Let $\kra$ be the bijection between $\S_n(132)$ and $\D_n$
described above. It follows from the above reasoning (and is also
shown in \cite{Eli02}) that $\kra$ maps fixed points to centered
tunnels, and excedances to right tunnels, i.e.,
$\fp(\pi)=\ct(\kra(\pi))$ and $\exc(\pi)=\rt(\kra(\pi))$.
Therefore, the left hand side of (\ref{fpexc132}) equals
$\sum_{D\in\D} x^{\ct(D)} v^{\rt(D)} z^{|D|}$. The result now is
obtained just applying Theorem~\ref{gftun} for $u=w=1$.
\end{proof}

Compare this expression~(\ref{fpexc132}) with equation (2) in
\cite{Eli02}. Note that Theorem~\ref{th:old} follows from these
two expressions.

As a further application, we can use the bijection $\bij$ to give
the following refinement of Corollary~\ref{cor}, which gives an
expression for the multivariate generating function for number of
fixed points, number of excedances, and number of descents in
132-avoiding permutations. An analogous result for 321-avoiding
permutations is given in \cite[section 3]{Eli02}.

\begin{theorem}\label{th:L} Let $$L(x,v,p,z):=1+\sum_{n\ge 1}
\sum_{\pi\in\S_n(132)}{x^{\fp(\pi)}v^{\exc(\pi)}}p^{\des(\pi)+1}z^n.$$
Then \bea\label{expr_L}
L(x,v,p,z)=\frac{2(1+xz(p-1))}{1+(1+v-2x)z-vz^2
(p-1)^2+\sqrt{f_1(v,z)}},\eea where $f_1(v,z)=
1-2(1+v)z+[(1-v)^2-2v(p-1)(p+3)]z^2-2v(1+v)(p-1)^2 z^3+v^2 (p-1)^4
z^4$.
\end{theorem}

\begin{proof}
We use again that $\kra$ maps fixed points to centered tunnels,
and excedances to right tunnels. It can easily be checked that it
also maps descents of the permutation to valleys of the
corresponding Dyck path. Clearly, the number of valleys of any
nonempty Dyck path equals the number of peaks minus one (in the
empty path both numbers are 0). Thus, $L$ can be expressed as
$$L(x,v,p,z)=\sum_{D\in\D} x^{\ct(D)} v^{\rt(D)} p^{\#\{\mbox{peaks of }D \}}
z^{|D|}.$$

By Theorem~\ref{corresp}, $\bij$ maps centered tunnels into hills
and right tunnels into even rises. Let us see what peaks are
mapped to by $\bij$. Given a peak $\u \d $ in $D\in\D$, $D$ can be
written as $D=A\u \d C$, where $A$ and $C$ are the parts of the
path before and after the peak respectively. This decomposition
corresponds to a tunnel of $D$ that goes from the beginning of the
$\u $ to the end of the $\d $. Assume first that the peak occurs
in the left half (i.e., $\len(A)<\len(C)$). When $D$ is read in
zigzag, the $\u $ opens a tunnel that is closed by the $\d $ two
steps later. This produces in $\bij(D)$ an up-step followed by a
down-step two positions ahead, that is, an occurrence of $\u \star
\d $ in the Dyck word of $\bij(D)$, where $\star$ stands for any
symbol (either a $\u $ or a $\d $).

If the peak occurs in the right half of $D$ (i.e.,
$\len(A)>\len(C)$), the reasoning is analogous, with the
difference that the $\d $ opens a tunnel that is closed by the $\u
$ two steps ahead. So, such a peak produces also an occurrence of
$\u \star \d $ in $\bij(D)$. Reciprocally, we claim that an
occurrence of $\u \star \d $ in $\bij(D)$ can come only from a
peak of $D$ either in the left or in the right half. Indeed, using
the notation from the procedure above describing the inverse of
$\bij$, an occurrence of $\u \star \d $ in $\bij(D)$ corresponds
to either an occurrence of $oc$ in the left half of $W$ or an
occurrence of $co$ in the right half of $W$. In both cases, the
algorithm given above will match these two letters $c$ and $o$
with each other, so they correspond to an occurrence of $\u \d $
in $D$.

If the peak occurs in the middle (i.e., $\len(A)=\len(C)$), then
by Lemma~\ref{lemma}, $\bij(A\u \d C)=\bij(AC)\u \d $, so it is
mapped to an occurrence of $\u \d $ at the end of $\bij(D)$.
Clearly we have such an occurrence only when $D$ has a peak in the
middle.

Thus, we have shown that peaks in $D$ are mapped by $\bij$ to
occurrences of $\u \star \d $ in $\bij(D)$ and occurrences of $\u
\d $ at the end of $\bij(D)$, or, equivalently, to occurrences of
$\u \star \d $ in $\bij(D)\d $ (here $\bij(D)\d $ is a Dyck path
followed by a down-step). Denote by $\oc(D)$ the number of
occurrences of $\u \star \d $ in $D\d $. This implies that $L$ can
be written as
$$L(x,v,p,z)=\sum_{D\in\D} x^{h(D)} v^{\er(D)} p^{\oc(D)} z^{|D|}.$$

We are left with a Dyck path enumeration problem, which is solved
in the following lemma. Let $J$ be defined in
Lemma~\ref{lemma_oc}. It is easy to see that we have $L(x,v,p,z) =
1 + J(x,1,p,1,v,p,z)$. Making use of (\ref{eq1}) and (\ref{eq2}),
it follows at once that $$L(x,v,p,z) = \frac{1 - xz + xpz}{1 - xz
- z K_1},$$ where $K_1$ is given by $$zK_1^2 - [1 - z- vz +
v(1-p)^2 z^2]K_1 + p^2 vz = 0.$$
From these equations we obtain (\ref{expr_L}).
\end{proof}

\begin{lemma}\label{lemma_oc}
Denote by $\ih(D)$ ($\fh(D)$) the number of initial (final) hills
in $D$ (obviously, their only possible values are 0 and 1).
Denote
by $\mu(D)$ the number of occurrences of $\u\star\d$
in $D$. Then the generating function $$J(x,t,s,u,v,q,z) := \sum
x^{h(D)} t^{\ih(D)} s^{\fh(D)} u^{\odr(D)} v^{\er(D)} q^{\mu(D)}
z^{|D|},$$ where the summation is over all nonempty Dyck paths, is
given by \bea\label{eq1} J(x,t,s,u,v,q,z) =\frac{uz[xts + (1-xu(1-t)(1-s)z)K]}{1
- xuz - uzK},\eea where $K$ is given by \bea\label{eq2} uz K^2 -
[1 - (u+v)z + uv(1-q)^2 z^2] K + q^2 vz = 0.\eea
\end{lemma}

\begin{proof} Every nonempty Dyck path has one of the following four
forms: $\u\d$, $\u\d B$, $\u A\d$, or $\u A\d B$, where $A$ and
$B$ are nonempty Dyck paths. The generating functions of these
four pairwise disjoint sets of Dyck paths are
\renewcommand{\theenumi}{\roman{enumi}}
\ben
\item $xtsuz$,
\item $xtuzJ(x,1,s,u,v,q,z)$,
\item $uzJ(1,q,q,v,u,q,z)$,
\item $uzJ(1,q,q,v,u,q,z)J(x,1,s,u,v,q,z)$,
\een respectively. Only the third factor in (iii) and (iv) needs an explanation:
the hills of $A$ are not hills in $\u A\d$; an initial (final)
hill in $A$ gives a $\u\u\d$ ($\u\d\d$) in $\u A\d$; an odd (even)
rise in $A$ becomes an even (odd) rise in $\u A\d$.

Consequently, the generating function $J$ satisfies the functional
equation \bea\label{eq3} J(x,t,s,u,v,q,z) = \hspace{10cm} \\ \nn
xtsuz + xtuz J(x,1,s,u,v,q,z) + uz J(1,q,q,v,u,q,z) + uz
J(1,q,q,v,u,q,z)J(x,1,s,u,v,q,z).\hspace{-8mm}\eea

From equation (\ref{eq3}) it is clear that, whether interested or
not in the statistics ``number of initial (final) hills", we had
to introduce them for the sake of the statistic marked by the
variable $q$. Also, without any additional effort we could use two
separate variables to mark the number of $\u\u\d$'s and the number
of $\u\d\d$'s, and obtain a slightly more general generating
function, although we do not need it here.

Denoting $H = J(x,1,s,u,v,q,z)$, $K = J(1,q,q,v,u,q,z)$, equation
(\ref{eq3}) becomes \bea\label{eq6} J = xtsuz + xtuzH + uzK +
uzHK.\eea Setting here $t=1$, we obtain \bea\label{eq7} H = xsuz +
xuzH + uzK + uzHK.\eea Solving (\ref{eq7}) for $H$ and introducing
it into (\ref{eq6}), we obtain (\ref{eq1}).

It remains to show that $K$ satisfies the quadratic equation
(\ref{eq2}). Setting $x=1$, $t=q$, $s=q$ in (\ref{eq6}), and
interchanging $u$ and $v$, we get \bea\label{eq8} K = q^2 vz + qvz
M + vz \wh{K} + vz M \wh{K},\eea where $M = J(1,1,q,v,u,q,z)$ and
$\wh{K}$ is $K$ with $u$ and $v$ interchanged, i.e. $\wh{K} =
J(1,q,q,u,v,q,z)$.

Now in (\ref{eq6}) we set $x=1$, $t=1$, $s=q$, and we interchange
$u$ and $v$, to get \bea\label{eq10} M = qvz + vz M + vz\wh{K} +
vzM\wh{K}.\eea Eliminating $M$ from (\ref{eq8}) and (\ref{eq10}),
we obtain \bea\label{eq11} vz(2qvz-q^2 vz+1-vz)\wh{K} + (vz-1)K +
vzK\wh{K} + q^2 vz = 0.\eea Finally, eliminating $\wh{K}$ from
(\ref{eq11}) and the equation obtained from (\ref{eq11}) by
interchanging $u$ and $v$, we obtain equation (\ref{eq2}). Note
that, as expected, $J$ is symmetric in the variables t and s and
linear in each of these two variables.
\end{proof}

From Theorem~\ref{th:L} one can see that the first terms of
$L(x,v,p,z)$ are $$1 + xpz + (vp^2 + x^2 p)z^2 + (v^2 p^2 + vp^2 +
xvp^3 + xvp^2 + x^3 p)z^3+\cdots,$$ corresponding to Dyck paths of
semilength at most 3 (or equivalently, to 321-avoiding permutations of length at most 3).


\begin{thebibliography}{99}

\bibitem{Deu99_} E. Deutsch, Dyck path enumeration, {\it Discrete Math.}
204 (1999), 167--202.

\bibitem{Eli02} S. Elizalde, Fixed points and excedances in restricted
permutations, preprint, arxiv:math.CO/0212221.

\bibitem{EliPak} S. Elizalde, I. Pak, Bijections for refined restricted
permutations, preprint, arxiv:math.CO/0212328.

\bibitem{FlSe98} P. Flajolet, R. Sedgewick, {\it Analytic combinatorics}
(book in preparation) (1998). (Individual chapters are available
as INRIA Research Reports 1888, 2026, 2376, 2956, 3162.)

\bibitem{Kra01} C. Krattenthaler, Permutations with restricted
patterns and Dyck paths, {\it Adv. Appl. Math.} 27 (2001),
510--530.

\bibitem{Rei02} A. Reifegerste, On the diagram of 132-avoiding
permutations, preprint, arxiv:math.CO/0208006.

\bibitem{RSZ02} A. Robertson, D. Saracino, D. Zeilberger, Refined
Restricted Permutations, preprint, arxiv:math.CO/0203033.

\bibitem{SeFl96}  R. Sedgewick, P. Flajolet,
{\it An Introduction to the Analysis of Algorithms},
Addison-Wesley, Reading, Massachusetts, 1996.

\bibitem{EC1} R. Stanley, {\it Enumerative Combinatorics}, vol. I,
Cambridge Univ. Press, Cambridge, 1997.

\end{thebibliography}
\end{document}